\documentclass[10pt,reqno]{article}
\usepackage{amsmath,amssymb,amsthm}
\usepackage{esint}
\usepackage{bm}
\usepackage{url}
\usepackage{subfigure}
\usepackage{graphicx,color}
\usepackage[sort&compress]{natbib}
\usepackage{algorithm}
\usepackage{algpseudocode}
\usepackage{algorithmicx}
\usepackage{hyperref}
\usepackage{float}
\usepackage{booktabs,multirow}
\usepackage{hhline}
\usepackage{setspace}

\usepackage{caption}
\usepackage{hyphsubst}
\usepackage{caption}

\usepackage{float}
\usepackage{filecontents}
\usepackage{hyperref}

\hypersetup{colorlinks,citecolor=blue,linkcolor=blue, urlcolor=blue}
\usepackage{pdflscape}
\usepackage{breakurl}

\newcommand{\RR}{\mathbb{R}}

\DeclareMathAlphabet{\itbf}{OML}{cmm}{b}{it}

\newtheorem{thm}{Theorem}[section]

\newtheorem{rem}[thm]{Remark}

\numberwithin{equation}{section}

\newcommand{\email}[1]{\protect\href{mailto:#1}{#1}}

\newcommand{\pathfigures}{Figures/}
\graphicspath{{\pathfigures}}

\def\bp{{{\itbf p}}}
\def\by{{{\itbf y}}}
\def\bx{{{\itbf x}}}

\def\bw{{{\itbf w}}}

\begin{document}

\title{Comments on \textquotedblleft Design of fractional-order variants of complex LMS and NLMS algorithms for adaptive channel equalization\textquotedblright
}

\author{
Shujaat Khan\footnotemark[2]
\and 
Abdul Wahab\footnotemark[1]\, \footnotemark[3]
\and
Imran Naseem\footnotemark[4]\, \footnotemark[5]
\and
Muhammad Moinuddin\footnotemark[6]\, \footnotemark[7]
}
\maketitle
\renewcommand{\thefootnote}{\fnsymbol{footnote}}
\footnotetext[1]{Corresponding Author. E-mail address:  \email{abdul.wahab@sns.nust.edu.pk}.}
\footnotetext[2]{Bio-Imaging, Signal Processing, and Learning Lab., Department of Bio and Brain Engineering, Korea Advanced Institute of Science and Technology, 291 Daehak-ro, Yuseong-gu, 34141, Daejeon, South Korea (\email{shujaat@kaist.ac.kr}).}
\footnotetext[3]{Department of Mathematics, School of Natural Sciences, National University of Sciences and Technology (NUST), Sector H-12, 44000, Islamabad, Pakistan (\email{abdul.wahab@sns.nust.edu.pk)}.}
\footnotetext[4]{School of Electrical, Electronic and Computer Engineering, The University of Western Australia, 35 Stirling Highway, Crawley, Western Australia 6009, Australia (\email{imran.naseem@uwa.edu.au)}.}
\footnotetext[5]{College of Engineering, Karachi Institute of Economics and Technology, Korangi Creek, 75190, Pakistan.}
\footnotetext[6]{Center of Excellence in Intelligent Engineering Systems (CEIES), King Abdulaziz University, Jeddah, Saudi Arabia (\email{mmsansari@kau.edu.sa)}.}
\footnotetext[7]{Electrical and Computer Engineering Department, King Abdulaziz University, Jeddah, Saudi Arabia.}
\renewcommand{\thefootnote}{\arabic{footnote}}

\begin{abstract}

The purpose of this note is to discuss some aspects of recently proposed fractional-order variants of complex least mean square (CLMS) and normalized least mean square (NLMS) algorithms in \textit{``Design of Fractional-order Variants of Complex LMS and Normalized LMS Algorithms for Adaptive Channel Equalization'' [Nonlinear Dyn. \textbf{88}(2), 839-858 (2017)]}. It is observed that these algorithms do not always converge whereas they have apparently no advantage over the CLMS and NLMS algorithms whenever they converge.  Our claims are based on analytical reasoning and are supported by numerical simulations.   

\end{abstract}



\section{Introduction}\label{Sec:Intro}
The least mean square (LMS) is a widely used algorithm in adaptive signal processing \cite{Haykin, LMSBook1}.  It has many variants to deal with assorted signals and environmental conditions \cite{CLMS, LMS11, QKLMS, NLMS, qLMS}.  
Recently,  some fractional-order variants of the CLMS and the NLMS  (henceforth acronymed as the FCLMS and the FNLMS resp.) are proposed in \cite{CFLMS} suggesting improved steady-state and convergence performances in an adaptive filtering framework.   A system input vector $\bx(k):=\left[x(k), x(k-1),\cdots, x(k-M+1)\right]^T$ (a zero-mean Gaussian process), with time index $k\in\mathbb{Z}$, is passed through an $M-$taps channel ($M\in\mathbb{N}$)  with finite impulse response filter $\mathbf{h}=[h_0, h_1,\cdots, h_{M-1}]^T$. To cancel out the effects of the channel, the noisy observation contaminated by a zero-mean white Gaussian noise $n(k)$, i.e.,  
\begin{align}
y(k) := \sum_{i=0}^M h_i x(k-i) +n(k) = h_0x(k)+\sum_{i=1}^M h_i x(k-1)+n(k),
\tag{1}
\end{align} 
is fed as an input to an $N$-taps ($\mathbb{N}\ni N>M$) equalizer filter with weight vector $\bw=[w_0,w_1,\cdots, w_{N-1}]^T\in\RR^N$ or $\mathbb{C}^N$. 
An input regression vector of the equalizer, 
$$
\by(k):=\left[y(k), y(k-1),\cdots, y(k-N+1)\right]^T,
$$ 
is formed with the design objective to adjust  $\bw$  so that the output of the equalizer resemble $x(k)$ with minimum  instantaneous error  $e(k):=x(k-\Delta)-\hat{x}(k-\Delta)$ during the data transmission. Here $ x(k-\Delta)=:d(k)$ is the desired response, $\hat{x}(k-\Delta):=\bw^H\by(k)$  is the estimated output of the equalizer and $\Delta\in\{0,1,\cdots\}$. Accordingly, the mean squared error (MSE) based objective function, 
\begin{align}
\min_{\bw} E\left[e^*(k)\left(d(k)-\bw^H\by(k)\right)\right],\tag{3}\label{eq:3}
\end{align}
is considered and solved using the FCLMS and the FNLMS. In above and hereinafter, the superposed $*$, $T$, and $H$ indicate complex conjugate, transpose, and Hermitian transpose, respectively.

In this note, we discuss some aspects of the design of the LMS variants, the FCLMS and FNLMS,  and also of the simulation setup in \cite{CFLMS}. It is observed that these algorithms have apparently no improvement over the CLMS and the NLMS. We present our main remarks in Section \ref{s:rem} and provide some supporting simulation results in Section \ref{s:sim}. We provide our final conclusions in Section \ref{Sec:Con}.

\section{Main Remarks}\label{s:rem}

In order to facilitate the ensuing discussion, the design premise, symbols, notations, and equation numbers used in this note are consistent with \cite{CFLMS}. The corrected equations are marked by a superposed asterisk.

\subsection{Objective Function for the FCLMS}\label{r:function}
The simplified objective function derived in \cite{CFLMS} from MSE \eqref{eq:3} is  
\begin{align}
J(\bw)=\sigma_\bx^2-2\bw^H\bp +\bw^H\mathbf{R}\bw, 
\tag{5}\label{eq:5}
\end{align}   
where $\bp$ is the cross-correlation between the input and the output, $\mathbf{R}$ is the auto-correlation $N\times N-$matrix of the output $\by$, and $\sigma_\bx^2=E[d(k)^2]$. Note that $\bx$, $\by$, $d(k)$, and $\bw$ are complex for the CLMS. 
\begin{proof}[\textbf{Claim 1}] The objective function $J(\bw)$ in \eqref{eq:5} considered for FCLMS is complex valued and does not correspond to the  MSE \eqref{eq:3}.
\end{proof}

The actual form of the objective function for the CLMS is well-known (see, e.g., \cite[Eq. 2.99]{LMSBook1}) and is given by
\begin{align}
J_{\rm cor}(\bw)=E[|d(k)|^2]-2\Re\left\{\bw^H\bp\right\} +\bw^H\mathbf{R}\bw.
\tag{5*}\label{eq:5*}
\end{align}  
Herein, $\Re$ and $\Im$ denote the real and the imaginary parts, respectively.  The functional $J_{\rm cor}(\bw)$ is convex (hyper-paraboloid) and possesses a unique minimum when $e(k) \to 0$ (see, e.g., \cite{Haykin}). 

Notice that the simplified expression \eqref{eq:5} can be obtained using the assumption 
$
E[x^*(k-\Delta)\by(k)] = E[x(k-\Delta)\by^*(k)].
$ 
However,  this is not possible for a complex system unless $ E[\Im\{x(k-\Delta)\by^*(k)\}]=0$ for all $k$, i.e., the cross-correlation  $\bp$ between $\bx$ and $\by$, is strictly real, which is a very strong assumption for complex system identification. On the other hand, $\sigma_x^2:=E[d(k)^2]$ as taken in \cite{CFLMS} is  complex valued. Indeed,
$$
\sigma_x^2 = \Re\left\{E\left[d(k)\right]\right\}^2-\Im\left\{E\left[d(k)\right]\right\}^2+2\iota \Re\left\{E\left[d(k)\right]\right\}\Im\left\{E\left[d(k)\right]\right\}.
$$
This is simply because of the fact that for any complex number $z=a+\iota b$, one has $z^2 = a^2-b^2+\iota ab$ where $\iota = \sqrt{-1}$.  Substituting this in the expression \eqref{eq:5} of $J(\bw)$ and rearranging the terms, one gets
\begin{align*}
J(\bw)=&E\left[|d(k)|^2\right]-2\Im\left\{E\left[d(k)\right]\right\}^2+2\iota \Re\left\{E\left[d(k)\right]\right\}\Im\left\{E\left[d(k)\right]\right\}
\\
& -2\Re\left\{\bw^H\mathbf{p}\right\}-2\iota \Im\left\{\bw^H\mathbf{p}\right\} +\bw^H\mathbf{R}\bw.
\end{align*}
On further simplification, one arrives at the form
\begin{align}
J(\bw)=&\Bigg(E\left[|d(k)|^2\right]  -2\Re\left\{\bw^H\mathbf{p}\right\}+\bw^H\mathbf{R}\bw\Bigg)
\nonumber
\\
 &+2\Bigg(\iota \Re\left\{E\left[d(k)\right]\right\}\Im\left\{E\left[d(k)\right]\right\}-\Im\left\{E\left[d(k)\right]\right\}^2-\iota\Im\left\{\bw^H\mathbf{p}\right\}\Bigg) 
\nonumber 
\\ 
=& J_{\rm cor}(\bw) + Z(k), 
\tag{a} \label{a}
\end{align}
where 
$$
Z(k):= 2\Bigg(\iota \Re\left\{E\left[d(k)\right]\right\}\Im\left\{E\left[d(k)\right]\right\}-\Im\left\{E\left[d(k)\right]\right\}^2-\iota\Im\left\{\bw^H\mathbf{p}\right\}\Bigg)\in\mathbb{C}.
$$

The remarks are as follows. 
\begin{rem}
\begin{enumerate}\label{Rem1}
\item[]
\item  Equation \eqref{a}  indicates that the objective functional $J(\bw)$ is complex valued for complex system  identification problems owing to $Z(k)$ in \eqref{a}. This is a contradiction to the fact that the MSE is a real-valued quadratic function. This justifies \textbf{Claim 1}.

 \item Even if we consider functional \eqref{a} as it is for the sake of argument then the first part $J_{\rm cor}$ corresponds to the standard MSE functional (see, Eq. \eqref{eq:5*}) and renders the optimal solution on minimization; (see, e.g., \cite{Haykin, LMSBook1}). The second component $Z(k)$ corresponds to a $k$-dependent complex error that corrupts standard MSE functional and makes the output complex. This impedes $J(\bw)$ to converge at the first place since it is unknown if $Z(k)\to 0$ as $k\to+\infty$. Even if it converges (i.e., $Z(k)\to 0$), the solution will converge to the optimal solution with high steady state error owing to the contribution of  $Z(k)$. 
 
\item In the light of the discussion above, any variant of the CLMS based on the objective function \eqref{eq:5} is expected to either diverge or completely fail. Hence, the proposed FCLMS may only work in the real cases when it is  simply a fractional-order variant of the LMS introduced in \cite{FLMS}. The performance of the fractional-order variants in \cite{FLMS} has already been debated in \cite{bershad2017comments, WS}, where it is observed that they have no advantage over the conventional LMS. 

\item As will be discussed later on (see Point $3$ in Remark \ref{Rem2}), there is a discrepancy between the pseudo-code implementation in \cite[Table 1]{CFLMS}  and the theoretical derivation of the FCLMS. Consequently, the aforementioned discussion is only relevant to the theoretical presentation of the algorithm which, in turn, differs from pseudo-code implementation.

\end{enumerate}
\end{rem}

\subsection{Fractional Calculus}\label{r:calculus}

The update rules for the FCLMS and FNLMS are, respectively, defined by 
\begin{align*}
&\bw(k+1):=\bw(k)+\frac{\mu_1}{2}\left[-\frac{\partial J[\bw(k)]}{\partial \bw}\right]+\frac{\mu_2}{2}\left[-\frac{\partial^\nu J[\bw(k)]}{\partial \bw^\nu}\right],
\\
&\bw(k+1):=\bw(k)+\frac{\mu_1}{2\|\by(k)\|^2}\left[-\frac{\partial J[\bw(k)]}{\partial \bw}\right]+\frac{\mu_2}{2\|\by(k)\|^2}\left[-\frac{\partial^\nu J[\bw(k)]}{\partial \bw^\nu}\right].
\end{align*}
Herein,  $\mu_1$ and $\mu_2$ are the controlling parameters for the integral and fractional updates, $0<\nu<1$ is the fractional-order,  and ${\partial^\nu}/{\partial \bw^\nu}$ is the fractional gradient with respect to $\bw$ defined in terms of the left Riemann-Liouville fractional derivative  
$ _0D_t^\nu$ as in \cite[Eq. (14)]{CFLMS}.  Using the formula (see, e.g.,  \cite{Kilbas})
\begin{align}
_0D_t^\nu t^z = \frac{\Gamma(z+1)}{\Gamma(z-\nu+1)}t^{t-\nu}, 
 \tag{32}\label{eq:32}
\end{align}
the fractional gradient term for the FCLMS is presented  in \cite{CFLMS}  as 
\begin{align}
\frac{\partial^\nu J(\bw(k))}{\partial \bw^\nu}= -\Gamma(2){\by^T}(k)e^*(k) \odot \frac{\bw^{1-\nu}_l(k)}{\Gamma(2-\nu)},
\label{eq:36}
\tag{36}
\end{align}
where $\Gamma$ represents the Euler's Gamma function 
$$
\Gamma (z) :=\int_{0}^\infty e^{-t}t^{z-1}dt, \quad z\in\mathbb{C},\,\, \Re\{z\}>0,
$$
the exponent on $\bw$ is component-wise, and $\odot$ is a component-wise vector multiplication. For the FNLMS, the same gradient term  is used without complex conjugate on $e(k)$.

\begin{proof}[\textbf{Claim 2}]
The expression \eqref{eq:36} is  invalid in both real and  complex cases.
\end{proof}
The expression \eqref{eq:36} is not justified in \cite{CFLMS}.  Note that $\bw$ is complex for the FCLMS. Therefore, the fractional gradient of the real-valued  (and thus, non-holomorphic) function $J_{\rm cor}:\mathbb{C}^N\to \RR$ with respect to $\bw\in\mathbb{C}^N$ should be calculated in the sense of Wirtinger calculus (see, for instance, \cite{Delgado}). The fractional gradient of the non-holomorphic function $J_{\rm cor}$ with respect to complex $\bw$ is undefined in mathematics to the best of our knowledge. 

Let us consider the real fractional gradient for the FNLMS and rigorously use fractional calculus for the Riemann-Liouville derivatives. Towards this end, we express $J(\bw)$ as 
\begin{align}
J(\bw(k))= \sigma_\bx^2-2\sum_{n=1}^N w_{n-1}(k)p_n(k)+\sum_{n,m=1}^N w_{n-1}(k)w_{m-1}(k) R_{nm}(k), 
\tag{b}
\label{A1}
\end{align} 
where $p_n$ and $R_{nm}$ are the components of the instantaneous cross-correlation vector $\bp$ and the autocorrelation matrix $\mathbf{R}$. 
In order to find the component fractional derivative $\partial^\nu J/ \partial w_{\ell-1}^\nu$, we re-arrange \eqref{A1} as 
\begin{align*}
J(\bw(k))=& \sigma_\bx^2-2\sum_{\substack{n=1\\ n\neq \ell}}^N w_{n-1}(k)p_n(k)-2w_{\ell-1}(k) p_\ell(k)
+\sum_{\substack{n,m=1\\ n\neq \ell,m\neq \ell}}^N w_{n-1}(k)w_{m-1}(k) R_{nm}(k)
\\
&+\sum_{\substack{n=1\\ n\neq \ell}} w_{n-1}(k)w_{\ell-1}(k) R_{n\ell}(k)
+\sum_{\substack{m=1\\ m\neq \ell}} w_{\ell-1}(k)w_{m-1}(k) R_{\ell m}(k) 
+  w_{\ell-1}^2(k)R_{\ell\ell}(k),
\end{align*} 
or equivalently 
\begin{align*}
J(\bw(k))=& \left[\sigma_\bx^2-2\sum_{n=1, n\neq \ell}^N w_{n-1}(k)p_n(k)+\sum_{\substack{n,m=1\\ n\neq \ell, m\neq \ell}}^N w_{n-1}(k)w_{m-1}(k) R_{nm}(k)\right]
\nonumber
\\
&
+2w_{\ell-1}(k)\left[\sum_{n=1, n\neq \ell}^N w_{n-1}(k) R_{n\ell}(k)-p_\ell(k)\right] +  w_{\ell-1}^2(k) R_{\ell\ell}(k).
\end{align*} 
Here, we have made use of the fact that $R_{nm}=R_{mn}$. Notice that the first term is constant with respect to $w_{\ell-1}$. Therefore, by the definition of the Riemann-Liouville derivative and invoking the rule \eqref{eq:32}, one arrives at 
\begin{align}
\frac{\partial^\nu J}{\partial w_{\ell-1}^\nu}
&= \frac{w_{\ell-1}^{-\nu}}{\Gamma(1-\nu)}\Bigg[\sigma_\bx^2-2\sum_{\substack{n=1\\ n\neq \ell}}^N w_{n-1}(k)p_n(k)+\sum_{\substack{n,m=1\\ n\neq \ell, m\neq \ell}}^N w_{n-1}(k)w_{m-1}(k) R_{nm}(k)\Bigg]
\nonumber
\\
&\quad
+\frac{2w^{1-\nu}_{\ell-1}(k)}{\Gamma(2-\nu)}\left[\sum_{\substack{n=1\\ n\neq \ell}}^N w_{n-1}(k) R_{n\ell}(k)-p_\ell(k)\right] 
+ \frac{2w^{2-\nu}_{\ell-1}(k)}{\Gamma(3-\nu)}R_{\ell\ell}(k).
\tag{c}\label{eq:c}
\end{align}
Thus, comparing \eqref{eq:c} with \eqref{eq:36}, one can easily see that a rigorus application of fractional calculus renders a significantly different result than rule \eqref{eq:36}. In fact, the chain rule for fractional derivatives is misleading and is much more involved than that for ordinary derivatives. We suggest interested readers to read the article by Tarasov \cite{Tarasov} or consult the monograph \cite{Kilbas}. The above discussion validates \textbf{Claim 2}.

\subsection{Schematic Issues}\label{r:design}

Let us ignore the mathematical observations discussed in Sections \ref{r:function}--\ref{r:calculus} for an instance. 
Precisely, assume  for the sake of argument that the update equations 
\begin{align}
\bw(k + 1) =& \bw(k) + \mu_1 e^*(k) \bw^T(k) + {\mu_2} G e^*(k) \by^T(k) \odot \left[\frac{\bw^{1-\nu}(k)}{\Gamma(2-\nu)}\right],\tag{37}\label{eq:37}
\\
\bw(k + 1) =& \bw(k) + {\mu_l} \frac{e(k) \by^T(k)}{\vert\vert\by(k)\vert\vert^2+\varepsilon} + \mu_f G \frac{e(k) \by^T(k)}{\vert\vert\by(k)\vert\vert^2+\varepsilon} \odot \left[\frac{\bw^{1-\nu}(k)}{\Gamma(2-\nu)}\right], \tag{45}\label{eq:45}
\end{align}
are constituted by an approximation.  Herein, $\mu_f$ is the fractional step size control parameter, $||\by(k)||$ is the norm of $\by(k)$,  $\varepsilon$ is a small parameter introduced to avoid zero denominators, and $G=\Gamma(3-\nu)/\Gamma(2-\nu)\Gamma(3)$. 

\begin{proof}[\textbf{Claim 3}]
The update rules defined by \eqref{eq:37} and \eqref{eq:45} render complex outputs for negative weight iterates.
\end{proof}

\begin{proof}[\textbf{Claim 4}]
There is a discrepancy between the theoretical analysis and the pseudo-code implementation for the FCLMS.
\end{proof}

Having defined the update rules as in \eqref{eq:37} and \eqref{eq:45}, the design of the fractional algorithms in \cite{CFLMS}  is similar to the fractional LMS algorithm in \cite{FLMS} (see, for instance, \cite[Eq. 15]{FLMS}). Accordingly, they inherit the shortcomings of the fractional LMS discussed in \cite{bershad2017comments}.  More specifically, we have the following remarks. 

\begin{rem}\label{Rem2}
\begin{enumerate}
\item[]
\item Rule \eqref{eq:45} indicates that the update $\bw(k+1)$ will become complex due to the presence of fractional powers if any component  $w_i(k)$ is negative. In that case, the FNLMS algorithm will not converge at all (see, e.g., Fig. \ref{FNLMS_neg}). If all the weights are positive, the FNLMS will either diverge (because the update may become negative during the process as LMS update takes a zigzag path towards an optimal solution; see \cite{Haykin})  or provide no improvement over the NLMS (see, e.g., Figs. \ref{Fig3}, and the discussions in articles \cite{bershad2017comments, WS}). This justifies \textbf{Claim 3} in the case of the FNLMS.

\item  Rule \eqref{eq:37} also substantiates that for negative values of $w_i$ the fractional term $\bw^{1-\nu}(k)$ will be complex.  However, the FCLMS may converge owing to the integral part of the update equation \eqref{eq:37} that corresponds to the CLMS algorithm.  Nevertheless, it will converge to a high steady-state residual error generated by the fractional term (see, e.g., Fig. \ref{Fig4}, and also Point 2 in Remark \ref{Rem1}).  This justifies \textbf{Claim 3} in the case of the FCLMS.

\item In contrast to all the theoretical discussion, the pseudo-code implementation in \cite[Table 1]{CFLMS} is done by just augmenting the corresponding update equation  of the CLMS  by the right hand side of Eq. \eqref{eq:36} and the rest of the code remains the same. This substantiates that the objective function $J_{\rm cor}$ is actually implemented in the  pseudo-code  \cite[Table 1]{CFLMS}.  More specifically, the theoretical analysis  is performed using $J(\bw)$ and the numerical implementation is performed using $J_{\rm cor}(\bw)$. However, it appears to be favorable in some situations where the convergence of the CLMS algorithm is not stymied by the fractional part of the update equation for the FCLMS. This can be observed in Fig. \ref{Fig4}. This validates our \textbf{Claim 4}.

\end{enumerate}
\end{rem}

\subsection{Simulation Setup}\label{r:simolution}

In adaptive signal processing, performance comparison between algorithms can be made on the bases of different criteria.  Three important measures of performance are: (i) convergence rate,  (ii) steady-state error, and (iii) computational complexity.  From equations \eqref{eq:37} and \eqref{eq:45}, it is observed that the fractional-order variants are computationally very expensive as they require additional number of computations. Therefore, we focus only on convergence and steady-state measures for our experiments. 

For a fair evaluation, the conventional algorithms and their proposed counterparts must be setup at either an equal convergence (for the steady-state performance comparison) or an equal steady-state (for the convergence performance).  Also, if one algorithm is supposed to perform better than the other then a  higher convergence rate at the cost of lower steady-state error must be shown. 

\begin{proof}[\textbf{Claim 5}]
The simulation results do not delineate the actual convergence trends of the fractional algorithms.
\end{proof}
We argue that with the simulation parameters used in \cite[Sect. 4.1]{CFLMS}, 
the LMS and the NLMS algorithms converge slowly (see \cite[Fig. 3-5]{CFLMS}).  Other observations are stated below.

\begin{rem}\label{Rem3}
\begin{enumerate}
\item[]

\item In  \cite[Fig. 3-5]{CFLMS}, the linear scale is used for $y-$axis which makes it difficult to compare the steady-state error.

\item The performance shown in \cite[Sect. 4.1]{CFLMS} for the FNLMS is not a representative convergence trend.  As per our numerical experiment (see Section \ref{SSec:4} and Figs. \ref{FNLMS_neg}--\ref{Fig3}), the FNLMS algorithm diverges with the simulation parameters used in  \cite[Sect. 4.1]{CFLMS}.
 
\item In a random desired weight vector scenario (\cite[Sect. 4.1]{CFLMS}), the information of the random distribution is not provided. Therefore, the results provided in \cite[Sect. 4.1]{CFLMS} are not reproducible.

\end{enumerate}
\end{rem}

These observations as well as \textbf{Claim 5} are further justified through numerical simulations in the next section.

\section{Simulations}\label{s:sim}
To evaluate the performance of the FNLMS and the FCLMS algorithms, we consider  the problem of system identification.  The algorithms are evaluated for two evaluation protocols:  (i) the system with negative desired weights under noisy environment with signal-to-noise ratio (SNR) of $10$dB;  (ii)  the system with all positive weights without any noise.  The source code of all the experiments are available online; see \cite{Source}.  

The NLMS, the CLMS, and their fractional order variants are configured to equal performance at $f=\nu=1$.  The performance of the FNLMS and the FCLMS is observed for $f=\nu=0.9$, $0.8$, $0.7$, $0.6$, $0.5$, and $0.4$.

For real inputs, we consider  a random signal of length $1000$ obtained from a zero-mean Gaussian distribution with variance $1$.  For the complex signal, a similar configuration is used but the signal is obtained from a circular complex Gaussian distribution instead of a real random source. 

The experiments are repeated for $1000$ independent rounds and mean results are reported.  For each independent round, the weights were initialized with zeros.  The performance of all the algorithms is evaluted on mean deviation (MD) which is the $\ell_1-$ norm of the difference between the sought and the obtained weights, i.e., 
\begin{equation*}
	\Delta \bw(k) = \frac{|\bw(k)- \hat{\bw}(k)|}{N},
\end{equation*}
where $\bw$ and $\hat{\bw}$ are the sought and approximated weight vectors at $k$th iteration, respectively. Here, $|\cdot|$ is the $\ell_1-$ norm and $N$ is the length of the filter vector.

\subsection{Performance Evaluation of the FNLMS}\label{SSec:4}

In \cite{CFLMS}, no example of the implementation of the FCLMS is provided. Moreover, examples of the FCLMS and FNLMS are not presented for negative sought weights. In this note, we choose two evaluation protocols, one with partly negative weights and one with positive weights (corresponding to the case study in \cite{CFLMS}). The details of these evaluation protocols are provided below.

\subsubsection{Evaluation Protocol 1}
\begin{figure}[!htb]
	\centering
	\includegraphics[width=0.45\textwidth]{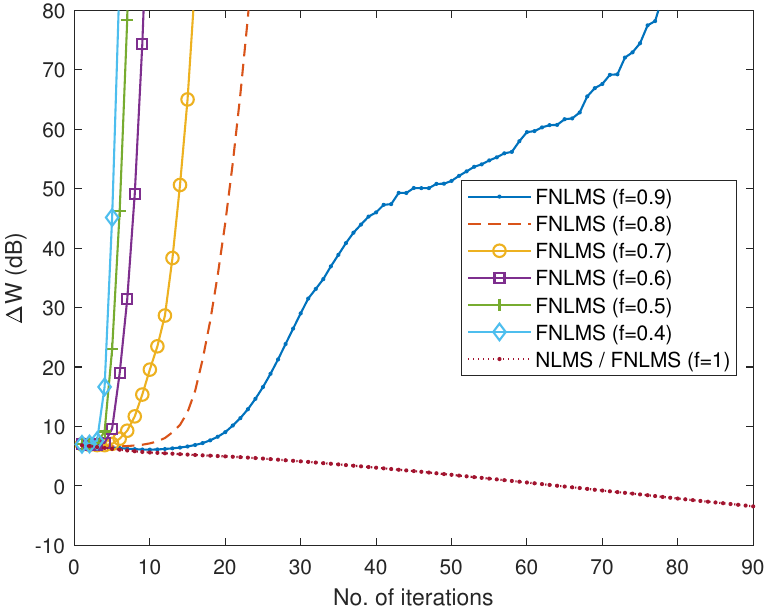}
	\caption{FNLMS: Learning curves for different values of fractional power ($f$) for negative weights in a noisy environment (SNR $=10$dB) with ${\mu_{\rm nlms}}=1$, $\beta=0.5$, and $\gamma=0.5$.}\label{FNLMS_neg}
\end{figure}

We consider a system with impluse response values 
$$
\bw=
\begin{bmatrix}
-10, &-9, &\cdots, &0, &\cdots, &9, &10
\end{bmatrix}^T\in\RR^{21}.
$$  
The step-size for the NLMS and FNLMS are set to ${\mu_{\rm nlms}}=1$, $\beta=0.5$, and $\gamma=0.5$ (see, \cite[Sect. 4]{CFLMS}). Figure \ref{FNLMS_neg}  shows the learning curves for the NLMS and the FNLMS.  We setup both algorithms at equal convergence rates, and compare  the steady-state performance of both algorithms.  In view of the simulations, it seems that the FNLMS algorithm completely fails to identify the system with negative weights for all the listed values of $\nu=f$ whereas the NLMS converges without any issue (see Fig. \ref{FNLMS_neg}).

\subsubsection{Evaluation Protocol 2}
\begin{figure}[!hbt]
\begin{center}
\subfigure[${\mu_{\rm nlms}=0.5}$]
{\includegraphics[width=.48\textwidth]{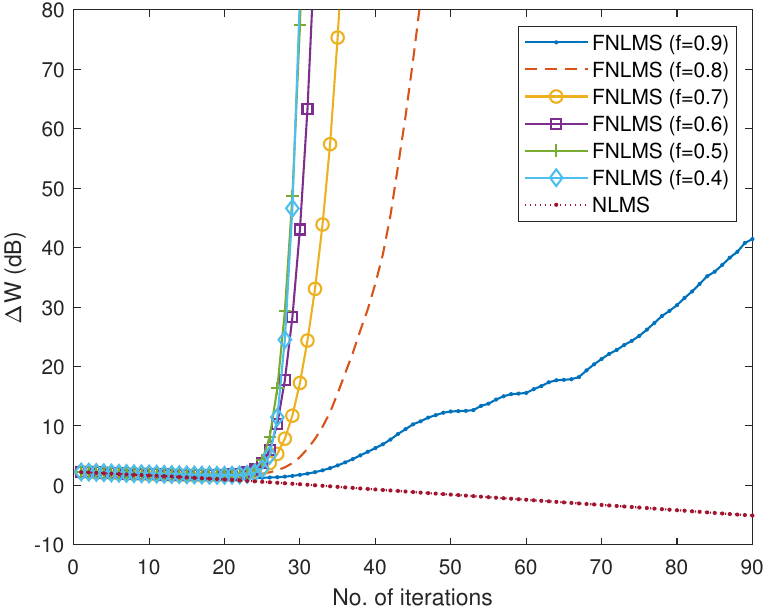}\label{FNLMS_pos1}}
\subfigure[${\mu_{\rm nlms}=1}$]
{\includegraphics[width=.48\textwidth]{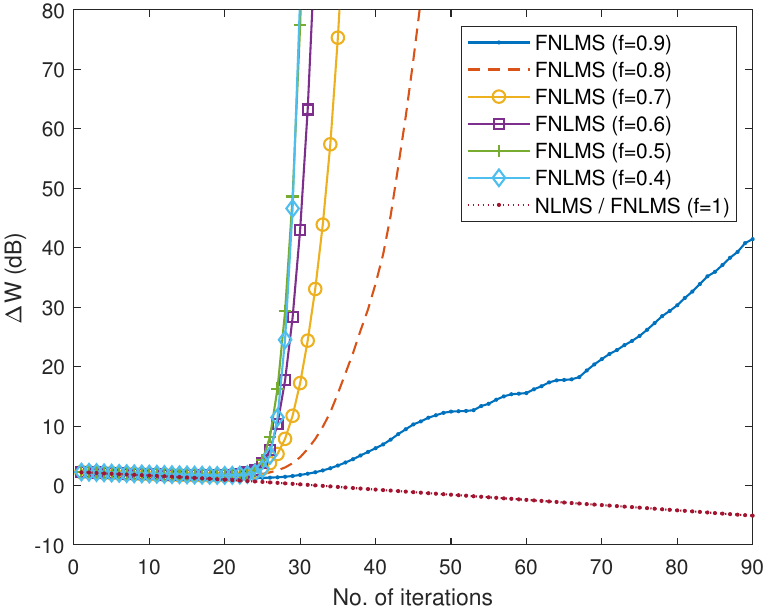}\label{FNLMS_pos2}} 
\caption{FNLMS: Learning curves for different values of fractional power ($f$) for positive weights in a noise-free case with $\beta=0.5$ and $\gamma=0.5$.}
\label{Fig3}
\end{center}
\end{figure}

In evaluation protocol 2, we choose the desired weight 
\begin{align*}
\bw=&
[1, 2, 2, 2, 1, 1, 2, 2, 3, 1, 1, 2, 2, 2, 1, 2, 1, 2, 2, 2, 1, 1, 2, 2, 2, 1, 1, 3, 2, 2]^T\in\RR^{30},
\end{align*}
and $\mu_{\rm nlms}=0.5$ for Fig. \ref{FNLMS_pos1} as in the case study performed in \cite[Sect. 4.1, Fig. 3]{CFLMS}.  Additionally, for Fig. \ref{FNLMS_pos2}, we choose $\mu_{\rm nlms}=1$ (when both algorithms are set up at an equal convergence rate).  The fractional step-sizes for both figures are set to be  $\beta=\gamma=0.5$. Both figures show the learning curves for the NLMS and the  FNLMS wherein the steady-state performances are compared.  From  Figs. \ref{FNLMS_pos1}-\ref{FNLMS_pos2}, it is observed that the FNLMS algorithm is diverging for all the listed values of $f=\nu$.

%

The results in Figs. \ref{FNLMS_neg} and \ref{Fig3} support \textbf{Claim 5}, i.e.,  the simulation results presented in \cite{CFLMS} do not show the representative convergence trends for the FNLMS algorithm.

\subsection{Performance Evaluation of the FCLMS}\label{Sec:4}

\begin{figure}[!hbt]
\begin{center}
\subfigure[Negative weights with noise (SNR$=10$dB)]
{\includegraphics[width=.48\textwidth]{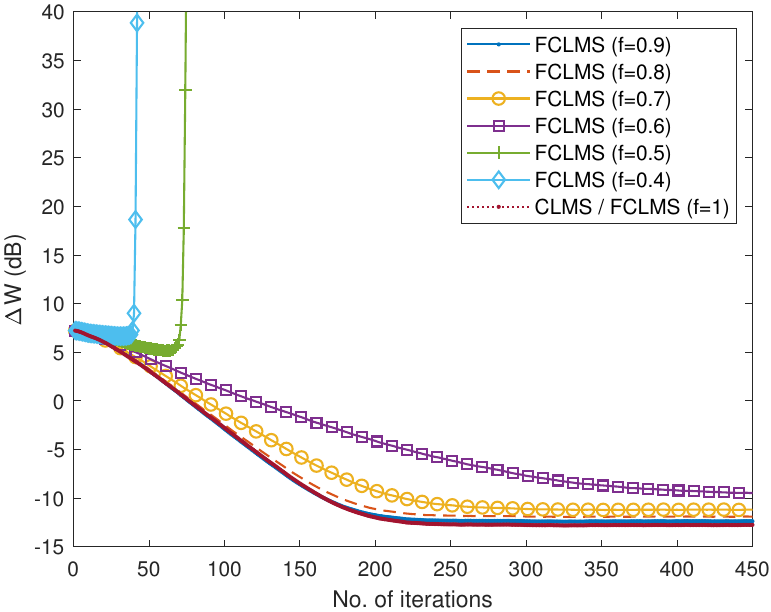}\label{FCLMS1}}
\subfigure[Positive weights without noise]
{\includegraphics[width=.48\textwidth]{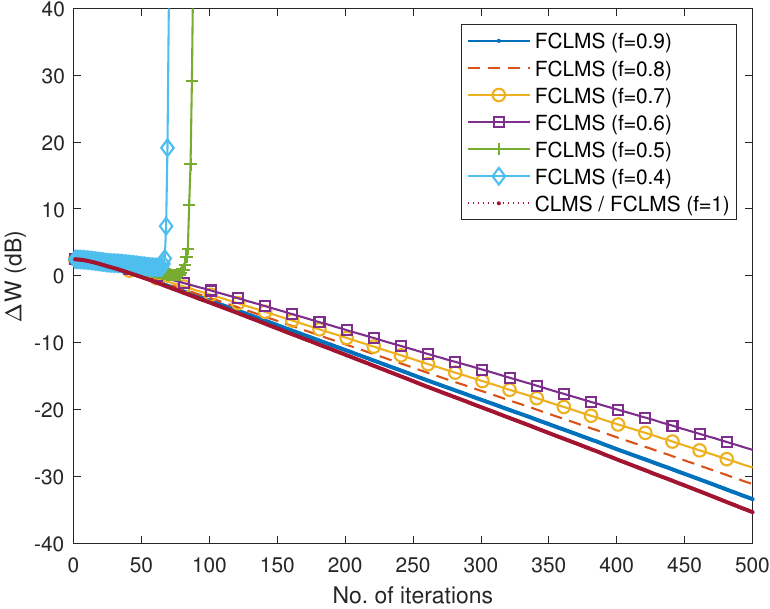}\label{FCLMS2}} 
\caption{FCLMS: Learning curves for different values of fractional power ($f$) with $\eta_{\rm clms}=0.04$, $\eta_{\rm fclms}=0.02$, and $\eta_f=0.02$.}
\label{Fig4}
\end{center}
\end{figure}

For the evaluation of the FCLMS, we  first consider a system with impulse response values 
$$
\bw=
\begin{bmatrix}
-10, &-9, &\cdots, &0, &\cdots, &9, &10
\end{bmatrix}^T\in\RR^{21}.
$$   
The step-size $\eta_{\rm clms}$ of CLMS is set to $0.04$, whereas the step-sizes for the FCLMS, $\eta_{\rm fclms}$ and $\eta_f$, are both set to $0.02$, respectively. Figure \ref{FCLMS1}, shows the learning curves for the CLMS and the FCLMS algorithms in a noisy environment with SNR$=10$dB.  We setup both algorithms at an equal convergence performance  and compare  the steady-state error. Fig. \ref{FCLMS1} indicates that the CLMS algorithm performs better than the FCLMS algorithm.  
Observe from these results that the fractional-term in the FCLMS has no advantage. Apparently, it is stymieing  the steady-state performance  of the integral part corresponding to the CLMS without even improving the convergence rate. 
In Fig. \ref{FCLMS2}, we perform another experiment without noise with desired weight vector
 \begin{align*}
\bw=&
[1, 2, 2, 2, 1, 1, 2, 2, 3, 1, 1, 2, 2, 2, 1, 2, 1, 2, 2, 2, 1, 1, 2, 2, 2, 1, 1, 3, 2, 2]^T\in\RR^{30},
\end{align*}
while keeping the rest of the setup identical. Similar conclusions as in the previous case hold. 
  
\section{Conclusion}\label{Sec:Con}
In this comment, we have analyzed the fractional-order variants of the complex least mean square (CLMS) and the normalized least mean square (NLMS) algorithms proposed in \cite{CFLMS}. We have discussed some aspects of the theoretical derivation, design, and simulation setup in \cite{CFLMS}. To be specific, the following points are highlighted:
\begin{enumerate}
\item The objective function considered for the FCLMS is complex valued and does not correspond to the MSE. 

\item The expression \cite[Eq. {\eqref{eq:36}}]{CFLMS} used for the fractional derivative is invalid in both real and complex cases.

\item The update rules of both the FCLMS and FNLMS render complex outputs for negative weight iterates. 

\item There is a discrepancy between the theoretical analysis and the pseudo-code implementation for the FCLMS. 

\item The simulation results do not delineate the actual convergence trends of the FCLMS and FNLMS algorithms.
\end{enumerate} 

The analysis performed in this note substantiates that the proposed algorithms either diverge or do not show any improvement in the performance in terms of convergence and steady-state error over the conventional algorithms.

\section{Conflict of Interest} 
The authors declare that they have no conflict of interest.

\bibliographystyle{plain}

\end{document}